\documentclass{amsart}
\usepackage{amsmath}
\usepackage{amscd}
\usepackage{amssymb}
\usepackage{amsfonts}
\newtheorem{theorem}{Theorem}[section]
\newtheorem{lemma}[theorem]{Lemma}

\newtheorem{proposition}[theorem]{Proposition}

\theoremstyle{definition}
\newtheorem{definition}[theorem]{Definition}

\theoremstyle{remark}
\newtheorem{remark}[theorem]{Remark}
\numberwithin{equation}{section}
\begin{document}

\title[The logarithmic entropy formula for the linear heat equation]
{The logarithmic entropy formula for the linear heat equation on Riemannian
manifolds}

\author{Jia-Yong Wu}
\address{Department of Mathematics, Shanghai Maritime University,
Haigang Avenue 1550, Shanghai 201306, P. R. China}

\email{jywu81@yahoo.com}

\date{\today}

\dedicatory{}


\begin{abstract}
In this paper we introduce a new logarithmic entropy functional for the
linear heat equation on complete Riemannian manifolds and prove that
it is monotone decreasing on complete Riemannian manifolds with
nonnegative Ricci curvature. Our results are simpler version, without
Ricci flow, of R.-G. Ye's recent result (arXiv: math.DG/0708.2008).
As an application, we apply the monotonicity of the logarithmic entropy
functional of heat kernels to characterize Euclidean space.
\end{abstract}
\maketitle
\section{Introduction}\label{sec1}
Given a compact $n$-dimensional Riemannian manifold $(M,g_0)$ without
boundary, the Ricci flow is the following evolution equation
\begin{equation}\label{flow}
\frac{\partial}{\partial t}g =-2Ric
\end{equation}
with the initial condition $g(x,0)=g_0(x)$, where $Ric$ denotes
the Ricci tensor of the metric $g(x,t)$. The Ricci flow equation was
introduced by R. Hamilton to approach the geometrization conjecture
in \cite{Hamilton}. Recently, studying various entropy functionals
along the Ricci flow is a very powerful tool for understanding of
Riemannian manifolds. A nice example is that G. Perelman
\cite{[Perelman]} introduced the following shrinking entropy functional
\[
\mathcal {W}(g(t), f(t),\tau):=\int_M \left[\tau\left(R+|\nabla
f|^2\right)+f-n\right](4\pi \tau)^{-n/2}e^{-f}d\mu,
\]
where $\tau>0$ and $d\tau/{dt}=-1$, $R$ and $d\mu$ denote the scalar
curvature and the volume form of the metric of $M$, respectively.
He proved that this entropy functional is nondecreasing along the Ricci flow
coupled to a backward heat-type equation. More precisely, if
$g(t)$ is a solution to the Ricci flow \eqref{flow} and the coupled
function $f(x,t)$ satisfies the evolution equation
\begin{equation}\label{hea1}
\frac{\partial f}{\partial t}=-\Delta f+|\nabla
f|^2-R+\frac{n}{2\tau},
\end{equation}
then Perelman proved that
\[
\frac{\partial \mathcal {W}}{\partial t}=2\tau\int_M\left|Ric
+\nabla^2f-\frac{g}{2\tau}\right|^2(4\pi\tau)^{-n/2}e^{-f}d\mu.
\]
The monotonicity of this entropy can be used to prove that shrinking
breathers must be shrinking gradient Ricci solitons. More
importantly, the monotonicity property is also fundamental in
proving Hamilton's little loop conjecture or what Perelman calls the
no local collapsing theorem (see \cite{[Perelman]} or \cite{[CCG]}).

\vspace{0.5em}

Another example is that in \cite{[FIN]}, M. Feldman, T. Ilmanen and L. Ni
tweaked some signs for Perelman's shrinking entropy $\mathcal {W}$,
and constructed a new entropy $\mathcal {W_{+}}$ corresponding to
expanding Ricci solitons, i.e.,
\[
\mathcal {W_{+}}(g(t), f_{+}(t),\sigma):=\int_M \left[\sigma\left(R+|\nabla
f_{+}|^2\right)-f_{+}+n\right](4\pi \sigma)^{-n/2}e^{-f_{+}}d\mu,
\]
where $\sigma>0$ and $d\sigma/{dt}=1$. They showed that this expanding entropy
$\mathcal {W_{+}}$ is also monotone nondecreasing on closed Riemannian manifolds
under the Ricci flow \eqref{flow} coupled to the backward heat-type
equation
\[
\frac{\partial f_{+}}{\partial t}=-\Delta f_{+}+|\nabla
f_{+}|^2-R-\frac{n}{2\sigma}
\]
and constant precisely on expanding Ricci solitons.

\vspace{0.5em}

Besides the above two entropies,  R.-G. Ye in \cite{[Ye]} also
introduced a new logarithmic entropy functional
\begin{equation}\label{logen1}
\mathcal{Y}_a(g(t),u(t),t):=-\int_M u^2 \log u^2 d\mu+\frac
n2\log\left[\int_M\left(|\nabla u|^2+\frac{R}{4}u^2\right)d\mu+a\right]+4at
\end{equation}
on an $n$-dimensional closed Riemannian manifold under the Ricci flow
\eqref{flow}, where $a$ is some constant. Here the coupled function
$u(t)$ satisfies
\begin{equation}\label{heateq1}
\frac{\partial u}{\partial t}=-\Delta u-\frac{|\nabla u|^2}{u}+\frac R2u
\end{equation}
such that $\int_M u^2 d\mu=1$. In other words, $u^2(t)$ solves to the
conjugate heat-type equation
\[
\frac{\partial u^2}{\partial t}=-\Delta u^2+Ru^2.
\]
satisfying $\int_M u^2 d\mu=1$. Under some suitable assumption, Ye can use
the monotonicity of Perelman's entropy to derive the monotonicity of
this logarithmic entropy functional along the Ricci flow \eqref{flow}
coupled to the heat-type equation \eqref{heateq1}. If we let
$u=(4\pi\tau)^{-n/4}e^{-f/2}$, then Perelman's entropy
functional can be rewritten as
\begin{equation}
\begin{aligned}\label{othform}
\mathcal {W}(g(t),u(t),\tau)=&-\int_Mu^2\log u^2d\mu
+4\tau\left[\int_M\left(|\nabla u|^2+\frac{R}{4} u^2\right)d\mu+a\right]\\
&-4a\tau-\frac n2\log (4\pi\tau)-n
\end{aligned}
\end{equation}
for an arbitrary constant $a$, and heat-type equation \eqref{hea1} becomes
equation \eqref{heateq1}. We would like to point out that Ye's
entropy \eqref{logen1} is often called \emph{logarithmic entropy}
because of the appearance of an additional logarithmic operation
compared with Perelman's entropy functional \eqref{othform}.

\vspace{0.5em}

The above mentioned entropy functionals are all considered with the
metric evolved by the Ricci flow. Below we recall an entropy
functional proposed by L. Ni \cite{[Ni1]} on a closed manifold
with a fixed metric. Let $(M,g)$ be an $n$-dimensional closed
Riemannian manifold with a fixed metric. Ni \cite{[Ni1]} (see
also chapter 16 in \cite{[CCG]}) considered the linear heat equation
\begin{equation}\label{heatequ}
\left(\frac{\partial}{\partial t}-\Delta\right)\tilde{u}=0
\end{equation}
and introduced the following entropy
\[
\mathcal {W}(f,\tau):=\int_M \left(\tau|\nabla
f|^{2}+f-n\right)(4\pi\tau)^{-n/2}e^{-f} d\mu,
\]
where $(f,\tau)$ satisfies
\[
\tilde{u}=\frac{e^{-f}}{(4\pi\tau)^{n/2}}\quad \mathrm{and}\quad
\int_M \frac{e^{-f}}{(4\pi\tau)^{n/2}}d\mu=1
\]
with $\tau>0$. By direct computation, Ni obtained the following
result.

\vspace{0.5em}

\noindent \textbf{Theorem A.} (L. Ni \cite{[Ni1]}) \emph{Let $(M,g)$ be
an $n$-dimensional closed Riemannian manifold. Assume that $\tilde{u}$
is a positive solution to the heat equation \eqref{heatequ} with
\[
\int_M \tilde{u} d\mu=1.
\]
Let the smooth function $f$ be defined as $\tilde{u}=(4\pi\tau)^{-n/2}e^{-f}$
and $\tau=\tau(t)$ with ${d\tau}/{dt}=1$. Then
\[
\frac{d \mathcal {W}}{dt}=-\int_M2\tau
\left(\left|\nabla^2f-\frac{g}{2\tau}\right|^{2}
+Ric(\nabla f,\nabla f)\right)\tilde{u}d\mu.
\]
In particular, if $M$ has nonnegative Ricci curvature, then
$\mathcal {W}(f,\tau)$ is monotone decreasing along the heat
equation \eqref{heatequ}.}

\vspace{0.5em}

We remark that if $\tilde{u}=H$ is the positive fundamental solution
of the heat equation \eqref{heatequ}, then Ni's entropy formula
on closed manifolds can be generalized to complete noncompact manifolds
with nonnegative Ricci curvature (see Lemma \ref{lem4.1} below).

In this paper, motivated by the work of Ye \cite{[Ye]}, we will
introduce a new logarithmic entropy functional for the linear heat
equation \eqref{heatequ} on a complete (possibly noncompact) manifold
under the static metric. This entropy functional is very similar to the
appearance of the logarithmic entropy functional along the Ricci flow
introduced by Ye \cite{[Ye]}. Following similar arguments to that of
Ye \cite{[Ye]}, we employ the property of Ni's entropy functional
and derive the monotonicity of our logarithmic entropy functional
for the heat equation \eqref{heatequ} as long as the Ricci curvature
of the manifold is nonnegative. As an application, on the noncompact
case, we apply the monotonicity of our logarithmic entropy functional
of heat kernels to characterize the Euclidean space. The main results
of this paper are Theorem \ref{T101} in Section \ref{sect2},
Theorem \ref{T101ap} in Section \ref{noncomp} and Theorem
\ref{app} in Section \ref{sect3b}.

The rest part of this paper is organized as follows. In Section \ref{sect2},
we first give some definitions of logarithmic entropies and then we state
Theorem \ref{T101}, which may be a natural generalization of Ye's
logarithmic entropy. After that, in Section \ref{sect3}, we give a
detailed proof of Theorem \ref{T101}. The proof mainly follows the
arguments of Ye's logarithmic entropy functional along Ricci flow
\cite{[Ye]}. The most difference is that our proof here makes use of
the monotonicity of Ni's entropy formula. In Section \ref{noncomp},
we generalize Theorem \ref{T101} to the complete noncompact setting
(see Theorem \ref{T101ap}). In Section \ref{sect3b}, we apply Theorem
\ref{T101ap} to give a characterization of Euclidean space. In
Section \ref{sect4}, we generalize Theorem \ref{T101} to the case
of the weighted heat equation.

\section{Monotonicity of the logarithmic entropy}\label{sect2}
Now we start to give several definitions, which are very
similar to Ye's definitions in the Ricci flow case.
\begin{definition}\label{definition1}
Let $(M,g)$ be an $n$-dimensional closed Riemannian manifold.
For any function $u\in W^{1,2}(M)$ satisfying $\int_M u^2 d\mu=1$,
we first define the \emph{logarithmic entropy} formula as follows:
\[
\mathcal{Y}_0(u):=-\int_M u^2\log u^2 d\mu+\frac
n2\log\left(\int_M|\nabla u|^2d\mu\right),
\]
where function $u$ satisfies $\int_M|\nabla u|^2d\mu>0$.

Then we define the \emph{logarithmic entropy} formula with
\emph{remainder} $a$ ($a>-\lambda$) as follows:
\[
\mathcal{Y}_a(u):=-\int_M u^2\log u^2 d\mu+\frac
n2\log\left(\int_M|\nabla u|^2d\mu+a\right),
\]
where $u\in W^{1,2}(M)$ and $\int_M u^2 d\mu=1$. Here
$\lambda$ denotes the first nonzero eigenvalue of the
Laplace operator.

In general the above constant $a$ in $\mathcal{Y}_a(u)$ is not arbitrary.
We add the restricted condition: $a>-\lambda$ to guarantee that
$\int_M|\nabla u|^2d\mu+a>0$.
\end{definition}

In the end, we will give the \emph{adjusted logarithmic entropy} formula
including another parameter $t$.
\begin{definition}\label{definition3}
Let $u$ and $a$ be the same as the above definition. We define
\[
\mathcal{Y}_a(u,t):=-\int_M u^2 \log u^2 d\mu+\frac
n2\log\left(\int_M|\nabla u|^2d\mu+a\right)-4at.
\]
\end{definition}
The form of this entropy is similar to \eqref{logen1}, with the most
difference is that here $\mathcal{Y}_a(u,t)$ is defined under the
static metric.

\vspace{0.5em}

For any closed Riemannian manifold with a fixed metric, let
$\lambda$ denote the first nonzero eigenvalue of the Laplace operator.
Let $u=u(x,t)$ be a smooth positive solution of the
following equation
\begin{equation}\label{evolu4}
\frac{\partial u}{\partial t}=\Delta u+\frac{|\nabla u|^2}{u}
\end{equation}
such that the normalization condition
\[
\int_M u^2 d\mu=1
\]
holds for all $t$. In fact, if we let $\tilde{u}=u^2$, then
$\tilde{u}$ satisfies heat equation \eqref{heatequ} with the
restraint condition
\[
\int_M \tilde{u}d\mu=1.
\]

\vspace{0.5em}

Now we state the monotonicity of the adjusted logarithmic entropy
formula on closed manifolds as follows.
\begin{theorem}\label{T101}
Let $M$ be an $n$-dimensional closed Riemannian manifold
with the nonnegative Ricci curvature. Then the adjusted
logarithmic entropy $\mathcal{Y}_a(u,t)$ ($a>-\lambda$)
is monotone decreasing along the heat-type equation
\eqref{evolu4} (namely, $u^2$ solves the heat equation
\eqref{heatequ}) with $\int_M u^2 d\mu=1$. More precisely,
\begin{equation*}
\begin{aligned}
-\frac{d\mathcal{Y}_a(u(t),t)}{dt}&\geq\frac{n}{4\omega}\int_M
\left[\left|-2\frac{\nabla^2u}{u}+2\frac{\nabla u\otimes\nabla
u}{u^2}-\frac{4\omega}{n}g\right|^{2}+4Ric\left(\frac{\nabla
u}{u},\frac{\nabla u}{u}\right)\right]u^2d\mu\\
&=\frac{n}{4\omega}\int_M
\left(\left|\bar{f}_{ij}-\frac{4\omega}{n}g_{ij}\right|^{2}
+R_{ij}\bar{f}_i\bar{f}_j\right)\frac{e^{-\bar{f}}}{(4\pi
t)^{n/2}}d\mu,
\end{aligned}
\end{equation*}
where $u=(4\pi t)^{-n/4}e^{-\bar{f}/2}$ and
$\omega:=\omega(u(x,t),a)=\int_M |\nabla u|^2 d\mu+a>0$.
\end{theorem}
We need to emphasize that in a general setting, the function $\bar{f}$ here
may be different from the function $f$ employed in the proof of Theorem
\ref{T101} given below. This function $\bar{f}$ is used for the purpose of
simplifying the expressions in the above formulas.

\begin{remark}
We can generalize this result to the case of noncompact manifolds, which will
be discussed in Section \ref{noncomp}.
\end{remark}

\begin{remark}
We can also generalize Theorem \ref{T101} to the case of the weighted
heat equation, which is treated in Section \ref{sect4}. For the
concepts of the weighted heat equation, the reader can refer to
\cite{[LD]}, \cite{[Wang]} or \cite{Wu10}.
\end{remark}

\section{Proof of Theorem \ref{T101}}\label{sect3}
Now we introduce Ni's entropy formula
\[
\mathcal {W}(f,\tau)=\int_M \left(\tau|\nabla
f|^{2}+f-n\right)\frac{e^{-f}}{(4\pi\tau)^{n/2}} d\mu,
\]
for the heat equation \eqref{heatequ} with $\int_M\tilde{u}d\mu=1$. Here
$\tilde{u}=u^2$ satisfying
\[
u=(4\pi\tau)^{-n/4}e^{-f/2},
\]
i.e.,
\begin{equation}\label{fangche3}
f=-\log u^2-\frac n2 \log\tau-\frac n2 \log(4\pi).
\end{equation}
So we can rewrite Ni's entropy as follows:
\begin{equation}
\begin{aligned}\label{wfunc}
\mathcal{W}(f,\tau)&=-\int_M u^2\log u^2 d\mu+(4\tau)\cdot
\left(\int_M|\nabla u|^2d\mu+a\right)-\frac n2\log(4\tau)\\
&\quad-4a\tau-\frac n2\log\pi-n,
\end{aligned}
\end{equation}
where $a$ is an arbitrary constant.

At first, we have the following useful fact.
\begin{lemma}\label{lemm1}
Assume $a>-\lambda$. Let $u\in W^{1,2}(M)$ with
\[
\int_M u^2d\mu=1.
\]
Then the minimum of the function
\[
h(s)=s\left(\int_M|\nabla u|^2d\mu+a\right)-\frac n2\log s
\]
for $s>0$ is given by
\[
\min h=\frac n2\log\left(\int_M|\nabla u|^2d\mu+a\right)+\frac
n2\left(1-\log\frac n2\right)
\]
and is achieved at the unique minimum point
\[
s=\frac n2\left(\int_M|\nabla u|^2d\mu+a\right)^{-1}.
\]
\end{lemma}
\begin{proof}
The proof can proceed essentially along the same line as in
\cite{[Ye]}. For completeness, we still present the proof in detail.
Let
\[
\omega=\int_M|\nabla u|^2d\mu+a.
\]
Then we have
\[
h(s)=\omega s-\frac n2\log s.
\]
Since $a>-\lambda$, we have $\omega>0$. So we know $h(s)\rightarrow
\infty$ as $s\rightarrow \infty$, and $h(s)\rightarrow \infty$ as
$s\rightarrow 0$. Therefore the function $h$ achieves its minimum at
the unique minimum point $s=\frac{n}{2\omega}$, since
\[
h'(s)=\omega-\frac{n}{2s}.
\]
Hence the minimum of $h$ is
\begin{equation*}
\begin{aligned}
h\left(\frac{n}{2\omega}\right)&=\frac n2-\frac n2\log
\left(\frac{n}{2\omega}\right)\\
&=\frac n2\log \omega+\frac n2\left(1-\log\frac n2\right).
\end{aligned}
\end{equation*}
This completes the proof of the lemma.
\end{proof}

Using this lemma we have a lower bound of Ni's entropy when $a>-\lambda$.
\begin{lemma}\label{lemm2}
Assume $a>-\lambda$. Then for each $\tau>0$, we have
\[
\mathcal{W}(f,\tau)\geq-\int_M u^2\log u^2 d\mu+\frac
n2\log\left(\int_M|\nabla u|^2d\mu+a\right)-4a\tau+b(n),
\]
where
\[
b(n):=-\frac n2\log\pi-\frac n2\left(1+\log\frac n2\right).
\]
Moreover, we have
\begin{equation}
\begin{aligned}\label{xiajie2}
\mathcal{W}\left(f,\frac{n}{8\omega(u,a)}\right)&=-\int_M u^2\log
u^2 d\mu+\frac n2\log\left(\int_M|\nabla
u|^2d\mu+a\right)\\
&\quad-\frac{na}{2\omega(u,a)}+b(n),
\end{aligned}
\end{equation}
where
\[
\omega(u,a):=\int_M|\nabla u|^2d\mu+a.
\]
\end{lemma}
\begin{proof}
This conclusion follows from the above lemma and \eqref{wfunc} immediately.
\end{proof}

\begin{remark}
Lemma \ref{lemm2} is still true for complete noncompact Riemannian
manifolds as long as $\mathcal{W}(f,\tau)$ is finite.
\end{remark}

\vspace{0.5em}

Now we can finish the proof of Theorem \ref{T101}.
\begin{proof}[Proof of Theorem \ref{T101}]
Let $u(x,t)$ be a smooth positive solution of equation
\eqref{evolu4}. Let $t_1\leq t_2$ and for $t\in[t_1,t_2]$,
we define for a given $\sigma>0$
\[
\tau=\tau(t)=t-t_1+\sigma.
\]
Assume that $f(x,t)$ is defined by \eqref{fangche3}, i.e.,
\[
u=(4\pi\tau)^{-n/4}e^{-f/2}.
\]
Obviously, $f(x,t)$ solves the
following equation
\[
\frac{\partial f}{\partial t}=\Delta f-|\nabla f|^2-\frac{n}{2\tau}.
\]
According to Ni's entropy monotonicity formula along the heat
equation, we have for $f=f(x,t)$ and $\tau=\tau(t)$
$({d\tau}/{dt}=1)$
\[
\frac{d \mathcal {W}}{dt}=-\int_M2\tau
\left(\left|f_{ij}-\frac{g_{ij}}{2\tau}\right|^{2}
+R_{ij}f_if_j\right)\tilde{u}d\mu
\]
on time interval $[t_1,t_2]$. Therefore
\[
\mathcal {W}(f(t_2),t_2-t_1+\sigma)-\mathcal
{W}(f(t_1),\sigma)=-2\int^{t_2}_{t_1}\int_M\tau
\left(\left|f_{ij}-\frac{g_{ij}}{2\tau}\right|^2
+R_{ij}f_if_j\right)\tilde{u}d\mu dt.
\]
If we choose
\[
\sigma=\frac{n}{8\omega(u(x,t_1),a)},
\]
then the above equality becomes
\begin{equation}
\begin{aligned}\label{nires3}
\mathcal{W}&(f(t_2),t_2-t_1+\sigma)+2\int^{t_2}_{t_1}\int_M\tau
\left(\left|f_{ij}-\frac{g_{ij}}{2\tau}\right|^{2}
+R_{ij}f_if_j\right)\tilde{u}d\mu dt\\
&=\mathcal {W}(f(t_1),\sigma)\\
&=-\int_M u^2\log u^2 d\mu\Big|_{t_1}+\frac
n2\log\left(\int_M|\nabla
u|^2d\mu+a\right)\Big|_{t_1}-4a\sigma+b(n),
\end{aligned}
\end{equation}
where we used \eqref{xiajie2}. On the other hand, by Lemma
\ref{lemm2}, we notice that
\begin{equation*}
\begin{aligned}
\mathcal{W}(f(t_2),t_2-t_1+\sigma)&\geq-\int_M u^2\log u^2
d\mu\Big|_{t_2}+\frac n2\log\left(\int_M|\nabla
u|^2d\mu+a\right)\Big|_{t_2}\\
&\,\,\,\,\,\,-4a(t_2-t_1+\sigma)+b(n).
\end{aligned}
\end{equation*}
Combining this with (\ref{nires3}) yields
\begin{equation*}
\begin{aligned}
-\int_M& u^2\log u^2 d\mu\Big|_{t_1}+\frac
n2\log\left(\int_M|\nabla u|^2d\mu+a\right)\Big|_{t_1}-4a\sigma+b(n)\\
&\geq2\int^{t_2}_{t_1}\int_M\tau
\left(\left|f_{ij}-\frac{g_{ij}}{2\tau}\right|^{2}
+R_{ij}f_if_j\right)\tilde{u}d\mu dt-\int_M u^2\log u^2 d\mu\Big|_{t_2}\\
&\quad+\frac n2\log\left(\int_M|\nabla
u|^2d\mu+a\right)\Big|_{t_2}-4a(t_2-t_1+\sigma)+b(n).
\end{aligned}
\end{equation*}
It follows that
\begin{equation}\label{entrequ}
\mathcal{Y}_a(u(t_1),t_1)\geq\mathcal{Y}_a(u(t_2),t_2)+2\int^{t_2}_{t_1}\int_M\tau
\left(\left|f_{ij}-\frac{g_{ij}}{2\tau}\right|^2
+R_{ij}f_if_j\right)\tilde{u}d\mu dt.
\end{equation}
Therefore, we obtain
\begin{equation*}
\begin{aligned}
-\frac{d\mathcal{Y}_a(u(t),t)}{d t}&\geq2\sigma\int_M
\left(\left|f_{ij}-\frac{g_{ij}}{2\sigma}\right|^2
+R_{ij}f_if_j\right)\frac{e^{-f}}{(4\pi\sigma)^{n/2}}d\mu\\
&=\frac{n}{4\omega}\int_M
\left[\left|-2\frac{\nabla^2u}{u}+2\frac{\nabla u\otimes\nabla
u}{u^2}-\frac{4\omega}{n}g\right|^{2}+4Ric\left(\frac{\nabla
u}{u},\frac{\nabla u}{u}\right)\right]u^2d\mu.
\end{aligned}
\end{equation*}
Hence the desired theorem follows.
\end{proof}

\section{On the complete noncompact case}\label{noncomp}
When $(M^n,g)$ is complete and noncompact, all of our above discussions
hold in Section \ref{sect3} as long as the integrations by parts make
sense and all the integrals involved are finite. In an analogous way
as the above argument, we can show that the monotonicity of the
adjusted logarithmic entropy $\mathcal{Y}_a(u,t)$ for the heat
kernel, still holds for any complete Riemannian manifold
with nonnegative Ricci curvature.

First, we recall the following entropy formula on complete noncompact
manifolds (see Corollary 16.17 and Theorem 16.27 in \cite{[CCG]})
\begin{lemma}\label{lem4.1}
Let $(M,g)$ be an $n$-dimensional complete noncompact Riemannian
manifold with nonnegative Ricci curvature. For
$\tilde{u}=u^2=(4\pi\tau)^{-n/2}e^{-f}=H$, the heat kernel of
heat equation \eqref{heatequ} satisfying $\int_M \tilde{u}d\mu=1$
and $\tau=\tau(t)$ with ${d\tau}/{dt}=1$, Ni's entropy
quantity $\mathcal {W}(f,\tau)$ is finite. Moreover,
\[
\frac{d \mathcal {W}}{dt}\leq-\int_M2\tau
\left|\nabla^2f-\frac{g}{2\tau}\right|^{2}u^2d\mu.
\]
\end{lemma}
The above lemma has been carefully proved in \cite{[CCG]}, which
involves estimates for the heat kernel and its first and second
derivatives. On the other hand, we have the following known facts
(see Corollaries 16.15 and 16.16 in \cite{[CCG]}), which are useful for
our discussion on the noncompact case.
\begin{lemma}\label{lem4.2}
Let $(M,g)$ be an $n$-dimensional complete noncompact Riemannian
manifold with nonnegative Ricci curvature. For
$\tilde{u}=u^2=(4\pi\tau)^{-n/2}e^{-f}=H$, the heat kernel of
heat equation \eqref{heatequ} satisfying $\int_M \tilde{u}d\mu=1$
and $\tau=\tau(t)$ with ${d\tau}/{dt}=1$, there exists a
constant $C(n)<\infty$ such that
\[
\int_Mf(x,y,\tau)H(x,y,\tau)d\mu\leq C(n)
\]
and
\begin{equation}\label{fjgj1}
\int_M|\nabla f|^2Hd\mu<\infty
\end{equation}
for any $x\in M$ and $\tau>0$.
\end{lemma}
Note that estimate \eqref{fjgj1} follows by Corollary 16.16 in \cite{[CCG]}
and Li-Yau's heat kernel upper bounds (Corollary 3.1 in \cite{[Li-Yau]}).
By Lemma \ref{lem4.2}, we easily have the following proposition.
\begin{proposition}\label{prop1}
Let $(M,g)$ be an $n$-dimensional complete noncompact Riemannian
manifold with nonnegative Ricci curvature. For
$\tilde{u}=u^2=(4\pi\tau)^{-n/2}e^{-f}=H$, the heat kernel of
heat equation \eqref{heatequ} satisfying $\int_M \tilde{u}d\mu=1$
and $\tau=\tau(t)$ with ${d\tau}/{dt}=1$, the adjusted logarithmic
entropy quantity $\mathcal{Y}_a(u,t)$ in Definition
\ref{definition3} (for any $x\in M$, $t\in \mathbb{R}$
and $a>-\lambda$) is finite.
\end{proposition}
\begin{proof}
Since
\[
-\log H=f+\frac n2\log (4\pi\tau)
\]
then we have
\begin{equation}\label{dejia1}
-\int_M u^2 \log u^2d\mu=\int_M f Hd\mu+\frac n2\log (4\pi\tau)
\end{equation}
and
\begin{equation}\label{dejia2}
\int_M|\nabla u|^2d\mu=\frac 14\int_M|\nabla f|^2Hd\mu.
\end{equation}
We also notice that
\[
\mathcal{Y}_a(u,t)=-\int_M u^2 \log u^2d\mu+\frac
n2\log\left(\int_M|\nabla u|^2d\mu+a\right)-4at.
\]
Hence our conclusion easily follows by \eqref{dejia1}, \eqref{dejia2}
and Lemma \ref{lem4.2}.
\end{proof}
From Proposition \ref{prop1}, we know that the entropy $\mathcal{Y}_a(u,t)$
is well-defined for $u^2=H$, being the heat kernel, on complete noncompact
Riemannian manifolds with nonnegative Ricci curvature.
Using Lemma \ref{lem4.1}, we can apply the same trick as in the proof
of Theorem \ref{T101} and obtain the following monotonicity of the
adjusted logarithmic entropy for the heat kernel on complete
(possibly noncompact) Riemannian manifolds.
\begin{theorem}\label{T101ap}
Let $(M,g)$ be an $n$-dimensional complete (possibly noncompact)
Riemannian manifold with nonnegative Ricci curvature. Then the
adjusted logarithmic entropy $\mathcal{Y}_a(u,t)$ ($a>-\lambda$)
with $\tilde{u}=u^2=H$, the heat kernel of heat equation
\eqref{heatequ} satisfying $\int_Mu^2d\mu=1$, is monotone
decreasing. More precisely,
\begin{equation}
\begin{aligned}\label{mainfor2}
-\frac{d\mathcal{Y}_a(u(t),t)}{dt}&\geq\frac{n}{4\omega}\int_M
\left|-2\frac{\nabla^2u}{u}+2\frac{\nabla u\otimes\nabla
u}{u^2}-\frac{4\omega}{n}g\right|^{2}u^2d\mu\\
&=\frac{n}{4\omega}\int_M
\left|\bar{f}_{ij}-\frac{4\omega}{n}g_{ij}\right|^{2}
\frac{e^{-\bar{f}}}{(4\pi t)^{n/2}}d\mu,
\end{aligned}
\end{equation}
where $u=(4\pi t)^{-n/4}e^{-\bar{f}/2}$ and
$\omega:=\omega(u(x,t),a)=\int_M |\nabla u|^2 d\mu+a>0$.
\end{theorem}
\begin{remark}
In general, $\bar{f}$ may be different from the function $f$
in Lemma \ref{lem4.1}. Here function $\bar{f}$ is used for the purpose of
simplifying the expressions in the above formulas.
\end{remark}

\section{Application}\label{sect3b}
In this section, we will use the monotonicity of the entropy functional
$\mathcal{Y}_a(u(t),t)$ of the heat kernel obtained in Section
\ref{noncomp} to characterize Euclidean space.
\begin{theorem}\label{app}
Let $(M^n,g)$ be a complete Riemannian manifold with nonnegative
Ricci curvature and constant $a>-\lambda$, where $\lambda$ is the first
nonzero eigenvalue of Laplace operator. If
$\mathcal{Y}_a(u(t_1),t_1)=\mathcal{Y}_a(u(t_2),t_2)$ for some
$t_1<t_2$, with $u^2$ being the heat kernel, then $(M^n,g)$ is
isometric to Euclidean space.
\end{theorem}
\begin{proof}
Since $\mathcal{Y}_a(u(t_1),t_1)=\mathcal{Y}_a(u(t_2),t_2)$ for some
$t_1<t_2$, using Theorem \ref{T101ap},
the monotonicity formula \eqref{mainfor2} (see also \eqref{entrequ}) implies
that
\[
\bar{f}_{ij}-\frac{4\omega}{n}g_{ij}\equiv0\Longleftrightarrow
f_{ij}-\frac{g_{ij}}{2(t-t_1+\sigma)}\equiv0
\]
on $(t_1,t_2)$, where
\[
\sigma=\frac{n}{8\omega(u(x,t_1),a)}\quad \mathrm{and}\quad
\omega(u(x,t_1),a)=\int_M |\nabla u|^2 d\mu+a>0.
\]
In other words,
\[
f_{ij}-\frac{g_{ij}}{2t}\equiv0
\]
for all $t\in(\sigma,t_2-t_1+\sigma)$. That is, for any such $t$,
$\varphi:=4tf$ satisfies
\begin{equation}\label{bian}
\nabla_i\nabla_j\varphi=2g_{ij}
\end{equation}
From this, we see that $\varphi$ attains its minimum at some point
$O\in M$. We claim that
\[
\varphi(x)-\varphi(O)=d^2(x,O).
\]
To prove this claim, let $x\in M$ be any point and consider a unit
speed minimal geodesic $\gamma:[0,d(x,O)]\to M$ joining $O$ to $x$.
We have
\[
\frac{d^2}{ds^2}\varphi(\gamma(s))=\nabla_i\nabla_j\varphi\dot{\gamma}^i\dot{\gamma}^j=2.
\]
Since $\nabla\varphi(O)=\vec{0}$, we have
$\varphi(\gamma(s))=\varphi(O)+d^2(O,\gamma(s))$ for all
$[0,d(x,O)]$ and hence the claim follows.

Therefore taking the trace of \eqref{bian} yields
\[
\Delta (d^2)\equiv 2n.
\]
This, together with the assumption that Ricci curvature is
nonnegative, implies $(M^n,g)$ is isometric to Euclidean space.
\end{proof}

\begin{remark}
The observation in proof of Theorem \ref{app} that
$f_{ij}-\frac{g_{ij}}{2t}\equiv0$ for all
$t\in(\sigma,t_2-t_1+\sigma)$ implies that $(M,g)$
is isometric to the standard Euclidean space has been
made by different authors. For example, see Corollary 1.3
in \cite{[Ni1]}, Proposition 2 in \cite{[PeWy]},
Theorem 3 in \cite{[PRS]} or Theorem 7.1 in \cite{[Nab]}.
\end{remark}
\section{Further remarks}\label{sect4}
While we only considered the linear heat equation in the previous
sections, the arguments in Sections \ref{sect2} and \ref{sect3} are
in fact valid for the weighted heat equation as well.

In order to make a clear statement of our result for the weighted heat
equation, we need to recall some basic facts about the $m$-dimensional
Bakry-\'{E}mery Ricci curvature (please see \cite{[Ba]}, \cite{[BE]},
\cite{[BE2]} and \cite{[LD]} for more details). Let $(M,g)$ be an
$n$-dimensional Riemannian manifold, and $h$ be a $C^2$ function.
We define a symmetric diffusion operator
\[
L:=\Delta-\nabla h\cdot\nabla,
\]
which is the infinitesimal generator of the Dirichlet form
\[
\mathcal{E}(\varphi_1,\varphi_2)=\int_M(\nabla \varphi_1,\nabla
\varphi_2)d\nu, \quad\forall \,\, \varphi_1, \varphi_2\in
C_0^{\infty}(M),
\]
where $\nu$ is an invariant measure of $L$ given by
$d\nu=e^{-h}d\mu$. It is well-known that $L$ is self-adjoint with
respect to the weighted measure $d\nu$. Given a smooth metric
measure space $(M,g,e^{-h}d\mu)$, the $\infty$-dimensional
Bakry-\'{E}mery Ricci curvature by
\[
Ric(L):=Ric+Hess(h),
\]
where $Hess$ denotes the Hessian of the metric $g$. We also define
the $m$-dimensional
Bakry-\'{E}mery Ricci curvature of the diffusion operator $L$ as
follows:
\[
Ric_{m,n}(L):= Ric(L)-\frac{\nabla h \otimes \nabla
h}{m-n},
\]
where $m:=\mathrm{dim}_{BE}(L)\geq n$ is called the Bakry-\'{E}mery
dimension of $L$, which is a constant and is not necessarily to be
an integer. In general the number $m$ is not equal to the manifold
dimension $n$, unless the operator $L$ is the Laplace operator. If
$m=\infty$, then $Ric_{m,n}(L)=Ric(L)$.

A remarkable feature of  $Ric_{m,n}(L)$ is that Laplacian comparison theorems
hold for $Ric_{m,n}(L)$ in the metric measure space $(M^m,g,e^{-h}d\mu)$
that look like the case of Ricci tensor in a $m$-dimensional manifold
\cite{[LD]} (see also \cite{Wu10}, \cite{[LD3]} and \cite{WeiWy}).
When $h$ is constant, $Ric(L)$ and $Ric_{m,n}(L)$ both recover
the ordinary Ricci curvature.

Given a smooth metric measure space $(M,g,e^{-h}d\mu)$, we consider
the weighted heat equation
\begin{equation}\label{weheat}
\frac{\partial \tilde{u}}{\partial t}-L\tilde{u}=0.
\end{equation}
In June 2006, X.-D. Li (see also \cite{[LD2]}) introduced the following entropy formula
\begin{equation}\label{waequ}
\mathcal {W}(f,\tau):=\int_M \left(\tau|\nabla
f|^{2}+f-m\right)(4\pi\tau)^{-m/2}e^{-f}d\nu,
\end{equation}
where $m$ is a finite constant satisfying $m\geq n$,
and $(f,\tau)$ satisfies
\[
\tilde{u}=\frac{e^{-f}}{(4\pi\tau)^{m/2}}\quad \mathrm{and}\quad
\int_M \frac{e^{-f}}{(4\pi\tau)^{m/2}}d\nu=1
\]
with $\tau>0$. By the direct calculation, he obtained the following
result.

\vspace{0.5em}

\noindent\textbf{Theorem B.} (X.-D. Li \cite{[LD2]}) \emph{Let
$M$ be an $n$-dimensional closed Riemannian manifold.
Assume that $\tilde{u}$ is a positive solution to the weighted heat equation
\eqref{weheat} with $\int_M \tilde{u} d\nu=1$. Let $f$ be defined by
\[
\tilde{u}=(4\pi\tau)^{-m/2}e^{-f}
\]
and $\tau=\tau(t)$ satisfies ${d\tau}/{dt}=1$. Then
\begin{equation*}
\begin{aligned}
\frac{d \mathcal {W}}{dt}&=-\int_M2\tau
\left[\left|\nabla^2f-\frac{g}{2\tau}\right|^2+Ric_{m,n}(L)(\nabla f,
\nabla f)\right]\tilde{u}d\nu\\
&\quad-\int_M\left(\sqrt{\frac{2\tau}{m-n}}\nabla h\cdot\nabla
f+\sqrt{\frac{m-n}{2\tau}}\right)^2\tilde{u}d\nu.
\end{aligned}
\end{equation*}
In particular, if $Ric_{m,n}(L)\geq 0$, then $\mathcal {W}(f,\tau)$
is monotone decreasing along the weighted heat equation
\eqref{weheat} with $\int_M \tilde{u} d\nu=1$.}

\begin{remark}
In fact, Li \cite{[LD2]} proved that the above result holds on complete
(noncompact) manifolds when $\tilde{u}=H$, being the heat kernel of
the weighted heat equation \eqref{weheat}. The above Theorem B also
appeared in \cite{[Wang]}.
\end{remark}

\begin{proof}
For the convenience of the reader, we give a sketch of the proof of
this theorem. If we let
\[
\phi=-\log \tilde{u} \quad \mathrm{and}\quad \omega=2L\phi-|\nabla \phi|^2,
\]
then
\[
\phi_t=L\phi-|\nabla \phi|^2.
\]
Now using the Bochner formula, we can derive that
\[
\left(\frac{\partial}{\partial t}-L\right)\omega=-2|\nabla^2\phi|^2
-2Ric(L)(\nabla\phi,\nabla\phi)-2\langle\nabla \omega,\nabla \phi\rangle.
\]
If we set $U:=\tau(2Lf-|\nabla f|^2)+f-m=\tau\omega+f-m$, then
\begin{equation*}
\begin{aligned}
\left(\frac{\partial}{\partial t}-L\right)U&=
-2\tau Ric_{m,n}(L)(\nabla f,\nabla f)-2\tau
\left|\nabla^2f-\frac{g}{2\tau}\right|^2\\
-&\left(\sqrt{\frac{2\tau}{m-n}}\nabla h\cdot\nabla
f+\sqrt{\frac{m-n}{2\tau}}\right)^2-2\langle\nabla U,\nabla f\rangle.
\end{aligned}
\end{equation*}
Hence
\begin{equation}
\begin{aligned}\label{shguan}
\left(\frac{\partial}{\partial t}-L\right)(U\tilde{u})&=
-2\tau\tilde{u}R_{m,n}(L)(\nabla f,\nabla f)
-2\tau\tilde{u}\left|\nabla^2f-\frac{g}{2\tau}\right|^2\\
&\quad-\tilde{u}\left(\sqrt{\frac{2\tau}{m-n}}\nabla h\cdot\nabla
f+\sqrt{\frac{m-n}{2\tau}}\right)^2.
\end{aligned}
\end{equation}
Note that
\[
\mathcal {W}(f,\tau)=\int_M(U\tilde{u})d\nu \quad \mathrm{and}\quad
\int_M|\nabla f|^2\tilde{u}d\nu=\int_MLf\tilde{u}d\nu.
\]
Therefore the result follows by integrating the equality \eqref{shguan}
with respect to the measure $d\nu$.
\end{proof}

Following the above definitions of Section \ref{sect2}, we give the
corresponding \emph{adjusted logarithmic entropy} formula for the diffusion
operator $L$.
\begin{definition}\label{definition4}
Let $(M,g,e^{-h}d\mu)$ be a metric measure
space, where $M$ is an $n$-dimensional closed Riemannian manifold.
Let $\kappa$ denote the first nonzero eigenvalue of the diffusion operator $L$.
For any function $u\in W^{1,2}(M)$ satisfying $\int_M u^2d\nu=1$,
if $a>-\kappa$, then we define a new logarithmic entropy with respect
to this metric measure space
\[
\mathcal{H}_a(u,t):=-\int_M u^2 \log u^2 d\nu+\frac
m2\log\left(\int_M|\nabla u|^2d\nu+a\right)-4at.
\]
\end{definition}

Parallel to  Lemma \ref{lemm2}, we have the following property.
\begin{lemma}\label{wlem}
Assume that $a>-\kappa$. Then for each $\tau>0$,
the entropy $\mathcal{W}$ defined by \eqref{waequ} satisfies
\[
\mathcal{W}(f,\tau)\geq-\int_M u^2\log u^2 d\nu+\frac
m2\log\left(\int_M|\nabla u|^2d\nu+a\right)-4a\tau+c(m),
\]
where
\[
c(m):=-\frac m2\log\pi-\frac m2\left(1+\log\frac m2\right).
\]
Moreover,
\begin{equation*}
\begin{aligned}
\mathcal{W}\left(f,\frac{m}{8\omega(u,a)}\right)&=-\int_M u^2\log
u^2 d\nu+\frac m2\log\left(\int_M|\nabla
u|^2d\nu+a\right)\\
&\quad-\frac{ma}{2\omega(u,a)}+c(m),
\end{aligned}
\end{equation*}
where $\omega(u,a):=\int_M|\nabla u|^2d\nu+a$.
\end{lemma}

For any metric measure space $(M,g,e^{-h}d\mu)$, let $u=u(x,t)$
be a smooth positive solution of the following equation
\begin{equation}\label{diffevolu4}
\frac{\partial u}{\partial t}=Lu+\frac{|\nabla u|^2}{u}
\end{equation}
such that the normalization condition $\int_M u^2d\nu=1$
holds for all $t$. Using Lemma \ref{wlem} and the
arguments of proving Theorem \ref{T101}, we can prove the following theorem.
\begin{theorem}\label{T102}
Let $(M,g,e^{-h}d\mu)$ be a metric measure
space, where $M$ is an $n$-dimensional closed Riemannian manifold.
If the $m$-dimensional Bakry-\'{E}mery Ricci curvature
is nonnegative and $a>-\kappa$, then the adjusted logarithmic entropy
$\mathcal{H}_a(u,t)$ is monotone decreasing along the heat-type
equation \eqref{diffevolu4} with $\int_M u^2d\nu=1$. More precisely,
\begin{equation*}
\begin{aligned}
-\frac{d\mathcal{H}_a}{d t}&\geq\frac{m}{4\omega}\int_M
\left[\left|-2\frac{\nabla^2u}{u}+2\frac{\nabla u\otimes\nabla
u}{u^2}-\frac{4\omega}{m}g\right|^{2}+4Ric_{m,n}(L)\left(\frac{\nabla
u}{u},\frac{\nabla u}{u}\right)\right]u^2d\nu\\
&\quad+\int_M\left(\sqrt{\frac{m}{\omega(m-n)}}\cdot\frac{\nabla
h\cdot\nabla u}{u}-\sqrt{\frac{4\omega(m-n)}{m}}\right)^2u^2d\nu\\
&=\frac{n}{4\omega}\int_M
\left(\left|\bar{f}_{ij}-\frac{4\omega}{m}g_{ij}\right|^2
+Ric_{m,n}(L)(\nabla\bar{f},\nabla\bar{f})\right)
\frac{e^{-\bar{f}}}{(4\pi t)^{m/2}}d\nu\\
&\quad+\int_M\left(\sqrt{\frac{m}{4\omega(m-n)}}\nabla
h\cdot\nabla
\bar{f}+\sqrt{\frac{4\omega(m-n)}{m}}\right)^2\frac{e^{-\bar{f}}}{(4\pi
t)^{m/2}}d\nu,
\end{aligned}
\end{equation*}
where $u=(4\pi t)^{-m/4}e^{-\bar{f}/2}$ and
$\omega:=\omega(u(x,t),a)=\int_M |\nabla u|^2 d\nu+a>0$.
\end{theorem}
\begin{proof}
The proof is very similar to that of Theorem \ref{T101}. Hence we
leave this proof to the interested reader. Notice that here we
should apply Theorem B instead of Theorem A and Lemma
\ref{wlem} to prove this theorem.
\end{proof}

\section*{Acknowledgments}
The author would like to thank Professor X.-D. Li for providing
papers \cite{[LD3]} and \cite{[LD2]}. He would also like to
express his gratitude to an anonymous referee for valuable
suggestions which led to the improvement of the paper. This work is
partially supported by NSFC (No. 11101267) and the Science
and Technology Program of Shanghai Maritime University
(No. 20120061).


\end{document}